\documentclass[11pt,
]{article}

\usepackage{amssymb}
\usepackage{amsfonts,amstext,amsmath,amsthm,
latexsym,mathrsfs,amsbsy,mparhack}
\usepackage[bbgreekl]{mathbbol}

\usepackage{graphicx}
\usepackage[svgnames]{xcolor}

\usepackage{marginnote}

\usepackage[russian,english,greek]{babel}
\languageattribute{greek}{polutoniko}

\usepackage{esvect}

\input{fsgd.sty}

\begin{document}

\selectlanguage{english}

\title{Definable selector for $\fd02$ sets modulo countable}

\author{Vladimir Kanovei\thanks
{
IITP RAS, 
Bolshoy Karetny, 19, b.1, Moscow 127051, Russia. 
Partial support of RFBR grant 17-01-00705 acknowledged. 
{\tt kanovei@googlemail.com}. 
}
\vyk{
\and
Vassily Lyubetsky\thanks
{IITP RAS, 
Bolshoy Karetny 19, b.1, Moscow 127051, Russia. 
Partial support of RFBR grant 18-29-13037 acknowledged. 
{\tt kanovei@googlemail.com}. }
}
}

\maketitle

\begin{abstract}
\noi 
A set is effectively chosen 
in every class of $\fd02$ sets modulo countable. 





\end{abstract}




Let $\cnt$ be the equivalence relation of equality 
modulo countable, that is, $X\cnt Y$ iff the symmetric difference 
$X\sd Y$ is (at most) countable. 
Does there exist an \rit{effective selector}, \ie, 
an effective choice of an element in 
each \dd\cnt equivalence class of sets of certain type? 
The answer depends on the type of sets considered. 
For instance, the question answers in the positive for 
the class of closed sets in Polish spaces by picking the only 
perfect set in each \dd\cnt equivalence class of closed sets. 
On the other hand, 
effective selectors for $\cnt$ 
do not exist 
in the domain of $\Fs$ sets, \eg, in the 
Solovay model 
(in which the axiom of choice AC holds and all ROD\snos
{ROD = real-ordinal definable, OD = ordinal-definable.} 
sets are LM and have the Baire property) 
by \cite[Theorem~5.5]{1}. 

Our goal here is to prove that $\Fs$ is the best possible 
for such a negative result. 

\bte
\lam{mt}
There exists a definable selector for\/ $\cnt$ in the domain of\/ 
$\fd02$ sets in Polish spaces. 
{\rm($\fd02$ = all sets simultaneously $\Fs$ and $\Gd$.)} 
\ete

\bpf[Theorem] 
We'll make use of the following lemma.

\ble
\lam{fdL}
If\/ $X$ is a countable\/ $\Gd$ set in a Polish space then 
the\/ {\ubf closure}\/ $\clo X$
is countable. 
Therefore if\/ $X\cnt Y$ are\/ $\fd02$ sets then\/ 
$\clo X\cnt\clo Y$. 
\ele
\bpf[Lemma]
Otherwise $X$ is a countable dense $\Gd$ set in an uncountable 
Polish space $\clo X$, which is not possible.
\epF{Lemma}

{\ubf Difference hierarchy.}
%
It is known 
(see \eg\ \cite[22.E]{kDST}) 
that every $\fd02$ set $A$ in a Polish space $\dX$ admits 
a representation in the form 
$A=\bigcup_{\et<\vt}(F_\et\bez H_\et)$, where 
$\vt<\omi$ and 
$F_0\qs H_0\qs F_1\qs H_1\qs\dots F_\et\qs H_\et\qs\dots$ 
is a decreasing sequence of closed sets in $\dX$, defined by 
induction so that $F_0=\dX$, 
$H_\et=\clo{F_\et\bez A}$, 
$F_{\et+1}=H_\et\cap \clo{F_\et\cap A}$, 
and the intersection on limit steps. 
The induction stops as soon as $F_\vt=\pu$.

The key idea of the proof of Theorem~\ref{mt} 
is to show that if $A\cnt B$ are $\fd02$ sets then 
the corresponding sequences of closed sets\pagebreak[0] 
$$
\left.
\bay{l}
F_0^A\qs H^A_0\qs F^A_1\qs H^A_1\qs\dots F^A_\et\qs H^A_\et\qs\dots
\\[1ex]
F_0^B\qs H^B_0\qs F^B_1\qs H^B_1\qs\dots F^B_\et\qs H^B_\et\qs\dots
\eay
\right\} \quad
(\et<\vt=\vt^A=\vt^B),\snos
{A shorter sequence is extended to the longer one by 
empty sets 
if necessary.}
$$ 
satisfying $A=\bigcup_{\et<\vt}(F^A_\et\bez H^A_\et)$ and 
$B=\bigcup_{\et<\vt}(F^B_\et\bez H^B_\et)$ as above, 
also satisfy 
\ben
\fenu
\itlb{*}
$F^A_\et\cnt F^B_\et$ \ and \ $H^A_\et\cnt H^B_\et$ \ --- \ 
for all $\et<\vt$. 
\een
It follows that the perfect kernels\snos
{$\pk X$, the \rit{perfect kernel}, 
is the largest perfect subset of a closed set $X$.} 
$\pk{F^A_\et}$, $\pk{F^B_\et}$ coincide: 
$\pk{F^A_\et}=\pk{F^B_\et}$, and 
$\pk{H^A_\et}=\pk{H^B_\et}$ as well.
Therefore the sets 
$\Phi(A)=\bigcup_{\et<\vt}(\pk{F^A_\et}\bez \pk{H^A_\et})$ 
and $\Phi(B)$ coincide 
(whenever $A\cnt B$ are $\fd02$ sets), 
and 
$A\cnt\Phi(A)$ holds 
for each $\fd02$ set $A$, 
so $\Phi$ is a selector required,
ending the proof of the theorem.

Thus it remains to prove \ref{*}. 
We argue by induction. 

We have $F^A_0=F^B_0=\dX$ (the underlying Polish space).

Suppose that $F^A_\et\cnt F^B_\et$; prove that 
$H^A_\et\cnt H^B_\et$. 
By definition, we have $H^A_\et=\clo{F^A_\et\bez A}$ and 
$H^B_\et=\clo{F^B_\et\bez B}$, where 
$(F^A_\et\bez A)\cnt (F^B_\et\bez B)$ 
(recall that $A\cnt B$ is assumed), 
hence $H^A_\et\cnt H^B_\et$ holds by Lemma~\ref{fdL}.

It's pretty similar to show that if 
$F^A_\et\cnt F^B_\et$ (and then 
$H^A_\et\cnt H^B_\et$ by the above) then  
$F^A_{\et+1}\cnt F^B_{\et+1}$ holds. 
This accomplishes the step $\et\to\et+1$.

Finally the limit step is rather obvious.
\epF{Theorem~\ref{mt}}

\bvo
Coming back to the mentioned result of \cite[Theorem~5.5]{1}, 
it is a challenging problem to prove that 
the equivalence relation $\cnt$ on $\Fs$ sets is not 
ROD-reducible to the equality of Bodel sets in the Solovay 
model.
\evo

\bre
As established in \cite{kl}, it is true in some models 
(including \eg\ Cohen and random extensions of $\rL$) that 
every OD and Borel set is OD-Borel (\ie, has an OD Borel code). 
In such a model, there is an effective choice of a set and 
its Borel code, 
by an OD function, in 
every \dd\cnt class of Borel sets containing an OD set. 
\ere

The author thanks Philipp Schlicht for useful comments.

\def\refname

\end{document}